\documentstyle[12pt]{article}
\newcommand{\const}{\mathop{\rm const}\limits}

\newcommand{\Ent}{\mathop{\rm Ent}\limits}

\begin{document}

\begin{center}

{\bf  THE RATE OF INCREASE FOR RECURSION \\

\vspace{3mm}

WITH QUADRATIC NON-LINEARITY } \par

\vspace{4mm}

 $ {\bf E.Ostrovsky^a, \ \ L.Sirota^b } $ \\

\vspace{4mm}

$ ^a $ Corresponding Author. Department of Mathematics and computer science, Bar-Ilan University, 84105, Ramat Gan, Israel.\\

E-mail: \ galo@list.ru \  eugostrovsky@list.ru\\

\vspace{4mm}

$ ^b $  Department of Mathematics and computer science. Bar-Ilan University,
84105, Ramat Gan, Israel.\\

E-mail: \  sirota3@bezeqint.net \\

\vspace{5mm}
                    {\bf Abstract.}\\

 \end{center}

 \vspace{4mm}

 We investigate in this short report the rate of increase of positive  numerical recursion with quadratic
non-linearity. More exactly, we intent to calculate the logarithmic index of its increasing. \par
 We present also the possible application in the theory of the Navier-Stokes equations. \par

\vspace{4mm}

\par

\vspace{4mm}

  \par

\vspace{4mm}

\section{Notations. Statement of problem. Conditions. Possible applications. }

\vspace{3mm}

 Let us consider the following numerical recurrence relation (dynamical system)

$$
D(n+1) = a + b \cdot D^2(n), \ a,b = \const > 0,  \ n = 0,1,2, \ldots, \ D(n) = D_{a,b}(n). \eqno(1.1)
$$

 We can and will suppose without loss of generality $ D(0) = 1  $ (initial condition). \par

 We impose in the sequel on the parameters $ a,b $ the following conditions:

 $$
b > 0;  \hspace{6mm} a\cdot b \ge 1/4. \eqno(1.2)
 $$
This conditions guarantee the monotonic increasing  of the sequence $ D(n): \\  D(n+1) \ge D(n). $\par

 Something similar occurred in the theory of elliptical curves and following in the coding theory
\cite{Hone1}. Another applications and investigations  of these equations are discussed in the book \cite{Agarwal1}. \par

\vspace{4mm}

{\it The authors are faced with this kind of equation by the investigation of numerical method for Navier-Stokes equation, see
\cite{Ostrovsky201}, \cite{Ostrovsky202}. } We describe in greater detail.  \par

\vspace{4mm}

 The mild solution $ u = u(x,t) $ of a Navier-Stokes equation in the whole space $  x \in R^d $ throughout
its lifetime $ t \in [0,T], \ 0 < T = \const \le \infty $ may be represented as a limit as $ n \to \infty, n = 0,1,2, \ldots $ the
following recursion:

 $$
 u_{n+1}(x,t) = u_0(x,t) + G[u_n, u_n](x,t), n=0,1,2,\ldots,
 $$
where $ u_0(x,t) $ is the solution of heat equation  with correspondent initial value and right-hand side
and $ G[u,v] $ is bilinear unbounded pseudo-differential operator, \cite{Kato1}, \cite{Kato2}.   See also the articles
\cite{Calderon1}, \cite{Cannone1}, \cite{Fabes1}, \cite{Giga1}, \cite{Hatami1}, \cite{Koch1}, \cite{Kozono1}, \cite{Serrin1},
\cite{Temam1} etc.  The second iteration is investigated  in \cite{German1}. \par

 Recall that the function $ u = u(x,t)  $  and hence the functions $ u_n(x,t), \ n = 0,1,2,\ldots $ are vector functions:

$$
u(x,t) = \vec{u}(x,t) = \{ u^{(i)}(x,t)  \}, \ i = 1,2,\ldots,d; \eqno(1.3)
$$
therefore the functional $ G[u,v] = G[\vec{u}, \vec{v}] = \vec{G}[\vec{u}, \vec{v}] $ may be interpreted as a tensor:

$$
\vec{G} = \{ g_{i,j}^m  \}, \
\vec{G}[\vec{u}, \vec{v}]_m = \sum_{i=1}^d \sum_{j=1}^d g_{i,j}^m u^{(i)} v^{(j)}, \ m = 1,2, \ldots,d. \eqno(1.4)
$$

 We denote by $ D(n) $ the amount of {\it independent} summands in the expression for the $ n^{th} $ iteration:

$$
u^{(i)}_n = \sum_{s_1=1}^{D(n)} \omega_{n, s_1}^{(i)}, \hspace{6mm} v^{(j)}_n = \sum_{s_2=1}^{D(n)} \kappa_{n, s_2}^{(j)}, \hspace{5mm}
\omega_{n,s_1}, \kappa_{n,s_2} =  \omega_{n,s_1}(x,t), \kappa_{n,s_2}(x,t), \eqno(1.5)
$$
 then

$$
\vec{G}[\vec{u_n}, \vec{v_n}]_m = \sum_{i=1}^d \sum_{j=1}^d g_{i,j}^m   \sum_{s_1=1}^{D(n)} \sum_{s_2=1}^{D(n)} \omega_{n,s_1}^{(i)} \kappa_{n,s_2}^{(j)}
u^{(i)}_n \ v^{(j)}_n, \ k = 1,2, \ldots,d. \eqno(1.6)
$$
 The last expression contains exactly in general case  $ d^2 \cdot D^2(n)  $ independent summands.\par

Obviously for all the values $  k  \  D(0) = 1 $ and

$$
D(n+1) = 1 +  d^2 \cdot D^2(n), \eqno(1.7)
$$
i.e. in this case $ a = 1, \ b = d. $  Since $  d \ge 1, $ the conditions (1.2) are satisfied.\par

\vspace{3mm}

{\bf  Our claim in this report is investigation of recurrence equation (1.1) under condition (1.2): obtaining of upper and lower
bounds and calculating the  asymptotic for the solution. } \par

\vspace{3mm}

\section{Main results: bilateral bounds and asymptotic behavior for solution.}

\vspace{3mm}

{\bf Theorem.}

$$
\forall k,l = 1,2,\ldots \ \Rightarrow
1 \le \frac{b \ D(k + l)}{ [ b \ D(l)]^{2^k}  } \le \left[ 1 + \frac{a}{b D^2(l)} \right]^{2^k - 1}; \eqno(2.1)
$$

$$
\forall k   \ge 1 \Rightarrow \lim_{l \to \infty} \frac{b \ D(k + l)}{ [ b \ D(l)]^{2^k}  } = 1. \eqno(2.2)
$$

\vspace{3mm}

{\bf Proof.}  {\it Lower bound:}

$$
D(k+1) \ge b \ D^2(k), \ k = l, l+1, l+2, \ldots; \ l = \const = 0,1,\ldots.
$$
 We deduce:

$$
D(l+1) \ge b \ D^2(l), \ D(l+2) \ge b^3 \ D^4(l), \ D(l+3) \ge b^7 \ D^8(l), \ldots
$$
 By induction:

$$
D(l+k) \ge b^{2^k - 1 } D^{ 2^k }(l). \eqno(2.3)
$$

\vspace{3mm}

{\it Upper bound.} We deduce denoting

$$
Q(l) = 1 + \frac{a}{b D^2(l)}
$$
and taking into account the monotonicity of the sequence $ D_{a,b}(n):  $ if $ k \ge l $ then

$$
D(k+1) = a + b D^2(k) = b D^2(k) \left( 1 + \frac{a}{b D^2(k)}  \right) \le
$$

$$
b D^2(k) \left( 1 + \frac{a}{b D^2(l)}   \right)  = b D^2(k) Q(l),
$$
and we find analogously

$$
D(k+l) \le b^{2^k -1} D^{ 2^k }(l) \left[  Q(l)   \right]^{ 2^k -1 }. \eqno(2.4)
$$

 The assertion (2.2) follows immediately from the bilateral estimates (2.1). \par

\vspace{3mm}

\section{Examples.}

\vspace{3mm}

 We intent to illustrate by building of some numerical examples the huge growth rate $ D(n) $ to infinity. \par

 \vspace{3mm}

{\bf 1.}  As regards to the Navier-Stokes equation in real case. \par

Here $  d = 3; $  i.e. $ D(n) = D_{1,9}(n); \ D(0) = 1, \ D(n+1) = 1 + 9 D^2(n):  $ \\

 $$
 D(0) = 1, \ D(1) = 10, \ D(2) = 901, \ D(3) = 811 \ 802, \ D(4) = 659 \ 022 \ 487 \ 205,
 $$

$$
D(5) = 434 \ 310 \ 638 \ 641 \ 864 \ 388 \ 712 \ 026, \ D(6) \approx 1.886257308 \cdot 10^{47},
$$

$$
D(7) \approx 3.5579666 \cdot 10^{94}.
$$

\vspace{3mm}

{\bf 2.} Let now $ a = b =1; $ i.e. $  D(n) = D_{1,1}(n);   $ then

$$
 D(0) = 1, \ D(1) = 2, \ D(2) = 5, \ D(3) = 26, \ D(4) = 677, \ D(5) = 458 \ 330,
$$

 $$
 \ D(6) = 210 \ 066  \ 388 \ 901, \ D(7) = 44 \ 127 \ 887 \ 745 \ 906 \ 175 \ 987 \ 802.
 $$

\vspace{3mm}

{\bf 3.}  For comparison:

$$
 D(n) \ge \tilde{D}(n) := 2^{(2^{n-1})};
 $$

\vspace{3mm}

 $$
 \tilde{D}(0) = 1, \  \tilde{D}(1) = 2, \  \tilde{D}(2) = 4, \  \tilde{D}(3) = 16,  \ \tilde{D}(4) = 256,  \tilde{D}(5) = 65 \ 536,
 $$

$$
\tilde{D}(6) = 4 \ 294 \ 967 \ 296, \  \tilde{D}(7)=  18 \ 446 \ 744 \ 073 \ 709 \ 551 \ 616.
$$

 The great difference between $  D_{1,1}(n) $ and $ \tilde{D}(n) $ show us the influence of free member "a" in  the
source equation (1.1). \par

\vspace{4mm}

\section{Concluding remarks.}

\vspace{3mm}

{\bf A.} At the same method may be used by investigation of the non-linear recursion

$$
D(n+1) = F(n, D(n))
$$
with monotonic increasing  power of non-linearity such that
$$
 C_1 z^{1 + \Delta} \le F(n,z)  \le C_2 z^{1 + \Delta}, \ z \ge 1, \ \Delta = \const > 0.
$$

\vspace{3mm}

{\bf B.} The vector analog of the equation (1.1) has a form

$$
D(n+1) = \vec{a} + D^2(n) \vec{b},
$$
where $ D(n)  $ is the square matrix $ m \times m $ and $ \dim  \vec{a} = \dim \vec{b} = m, \ m = 2,3,\ldots. $ \par

\vspace{3mm}

{\bf C.} Obviously, if in addition both the numbers $ a $ and  $ b $ are integer, then quite sequence $  \{   D(n) \}  $ is integer.
Therefore

$$
\forall k,l = 1,2,\ldots \ \Rightarrow b^{2^k - 1} \ D(l)^{2^k} \le
D(k + l) \le
$$

$$
\Ent \left\{ \left[ 1 + \frac{a}{b D^2(l)} \right]^{2^k - 1} \cdot b^{2^k - 1} \cdot D(l)^{2^k} \right\},
$$
where $ \Ent(z) $ denotes the integer part of the real number $ z, $ since $ D(n) $ is integer sequence. \par

\end{document}